\documentclass[12pt]{article}

\usepackage{amssymb,amsthm,amsmath}
\usepackage{url}

\def\eps{\varepsilon}
\def\ep{\varepsilon}

\def\E{\mathbb E}

\def\P{\mathbb P}

\def\NN{\mathcal N}

\newcommand{\R}{\mathbb{R}}

\newtheorem{lemma}{Lemma}[section]

\newtheorem{theorem}[lemma]{Theorem}

\newtheorem{cor}[lemma]{Corollary}

\newtheorem*{problem}{Problem}

\newcommand{\absconv}{\mathrm{absconv}}
\newcommand{\vol}{\mathrm{vol}}
\newcommand{\ptop}{V}
\newtheorem{fact}[lemma]{Fact}

\newcommand{\pr}{\noindent {\bf Proof: \, }}

\def\r{\right}

\def\no{\|\cdot\|}
\def\la{\langle}
\def\ra{\rangle}
\def\lam{\lambda}

\def\conv{\mathrm{conv}}
\def\absconv{\mathrm{absconv}}

\date{}

\title{On approximation by projections of polytopes with few facets}

\author{
Alexander E. Litvak${}^{1}$,
Mark Rudelson${}^{2}$,
Nicole  Tomczak-Jaegermann${}^{3}$}

\newcommand\address{\noindent\leavevmode

\noindent
Alexander E. Litvak, \\
Dept.~of Math.~and Stat.~Sciences,\\
University of Alberta, \\
Edmonton, Alberta, Canada, T6G 2G1.\\
\texttt{\small e-mail:  alexandr@math.ualberta.ca}

\medskip

\noindent
Mark Rudelson,\\
 Department of Mathematics,\\
 University of Michigan,\\
East Hall, 530 Church Street,\\
Ann Arbor, Michigan 48109\\
\texttt{\small e-mail:  rudelson@umich.edu}

\medskip

\noindent
Nicole  Tomczak-Jaegermann, \\
Dept.~of Math.~and Stat.~Sciences,\\
University of Alberta, \\
Edmonton, Alberta, Canada, T6G 2G1.\\
\texttt{\small e-mail:    nicole.tomczak@ualberta.ca}
}

\begin{document}

\maketitle

\footnotetext[1]{Research partially supported by  the
E.W.R. Steacie Memorial Fellowship.}
\footnotetext[2]{Research partially supported by NSF grant DMS 1161372.}
\footnotetext[3]{This author holds the Canada Research Chair in
  Geometric Analysis.}

\begin{abstract}
We provide an affirmative answer to a problem posed by Barvinok and
Veomett in \cite{BV}, showing that in general an $n$-dimensional convex body cannot 
be approximated by a projection of a section of a simplex of sub-exponential dimension. 
Moreover, we prove that for all $1\leq n\leq N$
there exists an $n$-dimensional convex body $B$ such that
for every $n$-dimensional convex body $K$ obtained as a projection of a
section of an $N$-dimensional simplex one has
$$
    d(B, K) \geq  c    \sqrt{  \frac{ n }{ \ln \frac{2 N \ln (2N)}{n}  }} ,
$$
where $d(\cdot, \cdot)$ denotes the Banach-Mazur distance
and $c$ is an absolute positive constant. The result is sharp up to
a logarithmic factor.

\bigskip

\noindent {\bf 2010 Subject Classification}:  Primary: 52A23, 52A27;  
\newline  Secondary: 52B55, 46B09.  \\ 


\end{abstract}

\section{Introduction}
\label{intro}

 One of the standard ways to describe a convex body in computational geometry
is the {\em membership oracle}.  The membership oracle of a body $K \subset \R^n$
is an algorithm, which, given a point $x \in \R^n$, outputs whether
$x \in K$, or $x \notin K$. If such oracle is constructed, and if the body $K$
has a relatively  well-conditioned position, meaning that
$r B_2^n \subset K \subset R B_2^n$ with $R/r \le n^C$, 
then one can construct
efficient probabilistic algorithms for estimating the volume of $K$, its inertia
ellipsoid, and other geometric characteristics 
(see e.g. \cite{KLS} and \cite{LVem}). Yet, constructing an efficient
membership oracle for a given convex body may be a hard problem \cite{BV}.
Because of this, it is important to know whether a convex body can be approximated
by another body, for which the membership oracle can be efficiently constructed.
One natural class of convex bodies for which the construction of the membership
oracle is efficient is the projections of a polytope with a few faces. Such polytopes
can be realized as projections of sections of a simplex in a dimension comparable to $n$.
This construction is discussed in details in \cite{BV}. In particular, the following
problem was posed (Problem 4.7.2 in \cite{BV}).

 \begin{problem}
    Let $K \subset \R^n$ be a symmetric convex body and let $P \subset \R^n$ be a
    projection of a polytope with $N$ facets, which approximates $K$ within a factor of
   $2$. Is it true that in the worst case the number $N$ should be at least exponential
    in $d$: $N \ge e^{cd}$ for some absolute constant $c>1$?
\end{problem}

Note that if $K=B_p^n$ is the unit ball of $\ell_p^n$, then this
approximation requires only proportional dimension. To see it recall
that a $(2n)$-dimensional simplex possesses a cubic section of
dimension $n$. Since a random projection of such a section is isomorphic
to an ellipsoid, we obtain an approximation of the Euclidean ball by a
projection of a section of a simplex in a dimension proportional to
$n$. Another deterministic construction of such an approximation was
found by Ben-Tal and Nemirovski \cite{BN}. A similar construction can
be used to approximate all balls $B_p^n$ for $2 \le p < \infty$. Since
the polar of a simplex is a simplex, one can also approximate the
balls $B_p^n$ for $1 \le p \le 2$. 
(Also, modifications of these constructions  give  explicit
symmetric ``conical subsets'' of the proportionally dimensional cube, 
whose linear projections can  arbitrarily close approximate the balls 
$B_p^n$ for $1 \le p <\infty$, see \cite{kt} for the details.)
Moreover,
even the existence of an $n$-dimensional convex
 body, which cannot be approximated by a projection of a section of a
 simplex $\Delta_N$ with $N$ proportional to $n$ has been an open
 problem.

 The main result of this paper provides an affirmative solution to 
the Barvinok problem above. Furthermore,
we prove a lower estimate for the minimal Banach--Mazur distance
between a certain convex symmetric body and a projection of a polytope
with $N$ facets. This estimate is optimal for all $N>n$ up to
logarithmic terms.

 \begin{theorem}\label{mainth}
Let $n\leq N$. There exists an $n$-dimensional convex symmetric body $B$, such that
for every $n$-dimensional convex body $K$ obtained as a projection of a
section of an $N$-dimensional simplex one has
$$
    d(B, K) \geq  c
   \sqrt{  \frac{ n }{ \ln \frac{2 N \ln (2N)}{n}  }} ,
$$
where $c$ is an absolute positive constant.
\end{theorem}

Let us note here that any projection of a section of
a simplex can be realized as a section of a projection of a simplex
and vice versa (see the next section). 
Thus, Theorem~\ref{mainth} holds for bodies $K$ obtained as a section of a
projection of a simplex as well.

To see that the estimate of Theorem~\ref{mainth} is close to optimal, recall that Barvinok 
proved in \cite{B1} that for every $N\geq 8 n$ and every symmetric convex body 
$B$ in $\R^n$ there exists a section $K$ of an $N$-dimensional simplex
such that
$$
  d(B, K) \leq C \max\left\{1, \sqrt{\frac{n}{\ln N}
        \cdot \ln \frac{n}{\ln N}}\right\}.
$$
Comparison of these two bounds shows that working with projections of 
sections of a simplex, as opposed to using sections alone, does not
significantly improve the approximation.  This is in stark contrast
with the situation described in the Quotient of a Subspace Theorem.
Recall that the Quotient of a Subspace Theorem of Milman (\cite{M},
see also \cite{MP} and \cite{R} for the non-symmetric case) states
that {\em given $\theta\in (0, 1)$ and an $n$-dimensional convex body
  $K$ there exists a projection of a section of $K$ whose dimension is
  greater than $\theta n$ and whose Banach-Mazur distance to the
  Euclidean ball of the corresponding dimension does not exceed
  $C(\theta)$ (moreover, $C(\theta)$ can be chosen such that
  $C(\theta)\to 1$ as $\theta \to 0^+$).}  On the other hand, it is
well-known by a volumetric argument (see Fact~\ref{volarg} below) that
any $n$-dimensional section of the $N$-dimensional cube (or simplex)
is at the distance at least $c\sqrt{n/\ln{(2N/n)}}$ from the
$n$-dimensional Euclidean ball. Thus, in the case of the cube (or
simplex) and proportional subspaces/projections, taking just sections
leads to $c\sqrt{n}$ distance to the Euclidean ball, while adding one
more operation -- taking a projection -- yields the distance bounded
by an absolute constant.

Our result also shows that Quotient of a Subspace Theorem cannot be
extended much beyond the Euclidean setting. Even if we start with the
simplest (in terms of complexity) convex body -- simplex -- we cannot
obtain an arbitrary convex set by taking a projection of a
section. Similar phenomena -- that many results of Asymptotic
Geometric Analysis cannot be extended much beyond the Euclidean
setting were discussed in \cite{LMT}.

It would be interesting to characterize the class of all
$n$-dimensional convex bodies, that can be realized (up to a
Banach-Mazur distance less than or equal to 2, say) as a projection of
a section of an $N$-dimensional simplex for $N=O(n)$. As we mentioned
above any $B_p^n$ is in this class, clearly any polytope with 
 $O(n)$ vertices or faces is in this class as well. In a related
direction we conjecture that there is no convex body $K$ such that an
arbitrary body can be obtained (up to Banach-Mazur distance bounded by
a constant)
from $K$ by taking a projection of a section.

Finally we would like to mention that many aspects of computational complexity of convex bodies
were discussed in \cite{S3}.

The paper is organized as follows.  In the next section we introduce
notation and auxiliary results, that will be used latter. 
We also describe a class of random polytopes crucial for our
construction in which we will find our example. 
We model these
polytopes on random polytopes introduced by Gluskin in \cite{G2}. 
In Section~\ref{main} we prove the main theorem, Theorem~\ref{mainth}.
The proof of this theorem uses Theorem~\ref{mt}, which states that
with high probability two Gluskin's polytopes are on large Banach-Mazur
distance to each other.  The last section is devoted to the proof of
Theorem~\ref{mt}.

\bigskip

\noindent
{\bf Acknowledgment.} The second author is grateful to Alexander
Barvinok for many helpful discussions.

\section{Notation and Preliminaries}
\label{notat}

By $|\cdot|$ and $\la \cdot , \cdot \ra$ we denote the canonical
Euclidean norm and the canonical inner product on $\R ^d$.  $B_2^d$
and $S^{d-1}$ stand for the Euclidean unit ball and the unit sphere,
respectively;  the standard basis of $\R ^d$ is denoted  by $e_1,
\ldots, e_d$.

As usual, $\| \cdot \| _p$, $1\leq p \leq \infty$, denotes the
$\ell _p$-norm, i.e. for every $x=(x_i)_{i=1}^d \in\R^d$
$$
 \|x\| _p = \left(
                   \sum _{i= 1}^d |x_i|^p
               \right) ^{1/p} \,
                                     \mbox{ for } \ p < \infty, \, \, \, \,
      \quad \, \, \, \, \|x\| _{\infty } = \sup _{i\leq d} |x_i|,
$$
and $\ell _p^d = (\R ^d, \|\cdot \|_p)$.
The unit ball of $\ell _p^d$ is denoted  by $B_p^d$.

Recall that
$\lceil x\rceil$ denotes the smallest integer which is not less than $x$.

By a convex body we mean a compact set with a non-empty interior.  For
a convex body $K\subset \R^d$ with $0$ in its interior, the Minkowski
functional of $K$ is
$$
    {\| x\|}_{K}=\inf \{\lambda >0 \ |\ x\in\lambda K\},
$$
i.e. it is the homogeneous convex functional, whose unit ball is
$K$. The polar of $K$ is
$$
 K^{\circ} = \{x \ | \ \la  x, y \ra \leq 1 \ \ \mbox{ for all } \ y\in K \}.
$$
Note that if $K$ is symmetric, then $K^{\circ}$ is the unit ball of the space dual to $(\R^d, \|\cdot\|_K)$.

It is well known that for any convex body  $K\subset \R^d$ there exists a point
$a\in K$ such that
\begin{equation}\label{center}
  -(K-a) \subset d (K-a).
\end{equation}
For example the center of the maximal volume ellipsoid contained in $K$
satisfies this (\cite{J}, see also \cite{Ba}).

Given a subset $K\subset \R^d$ the convex hull and the absolute convex hull 
of $K$ are denoted by $\conv (K)$ and $\absconv (K)=\conv (K\cup -K)$ respectively.  
The volume of $K$ is denoted by $\vol\ (K)$. 
A position of $K$ is a non-degenerate affine image of $K$.

For two convex bodies $K_1$ and $K_2$ in $\R^d$ the Banach-Mazur
distance between them is defined as
$$
  d(K_1, K_2) = \inf \{ \lam >0 \, \mid \, K_1 -a \subset
  T(K_2-b)\subset \lam ( K_1 -a ) \},
$$
where infimum is taken over all non-degenerate linear operators $T:
\R^d \to \R^d$ and all $a, b \in \R^d$.  Note that if $K_1=-K_1$ and
$K_2=-K_2$ then $a, b$ can be taken equal to $0$.  The distance 
$d(\cdot, \cdot)$ satisfies the multiplicative triangle inequality, i.e.  $d(K_1, K_2)\leq
d(K_1, K_3) d(K_3, K_2)$.

\smallskip

We  fix the following notation.
$$
    S:=S(N)= \left\{x =\{x_i\} _{i=1}^{N+1}\in \R^{N+1}\, \mid \,
      x_i\geq 0, i\leq N+1\r\},
$$
$$
    H:=H(N)= \left\{x =\{x_i\} _{i=1}^{N+1}\in \R^{N+1}\, \mid \, \sum
      _{i=1}^{N+1} x_i=1\r\},
$$
and
$$
    \Delta =\Delta_N := S\cap H.
$$
Note that $\Delta = \conv \{e_i\} _{i=1}^{N+1}$ is an $N$-dimensional
regular simplex.

As we mentioned in the introduction,  any projection of a section of
a simplex can be realized as a section of a projection of a simplex
and vice versa. Indeed, let $E$ be a linear, and let $F$ be an affine subspace of $\R^{N+1}$. 
Consider the body $K=P_E \Delta_{N+1} \cap F$. Without loss of generality, we 
may assume that $F \subset E$.  In this case $K=P_F (\Delta_{N+1} \cap \tilde{E})$, 
where $\tilde{E}=P_F^{-1}E=E \oplus F^\perp$.

Recall that a set $F\subset \R^{N+1}$ is an affine subspace if
there exists $b \in \R^{N+1}$ such that $F-b$ is a linear subspace of
$\R^{N+1}$. 
Given a set $K \subset \R^{N+1}$ and an affine subspace $F\subset
\R^{N+1}$   the section of $K$ by $F$ is denoted by
$$
    K^F = K\cap F.
$$
In particular,
$$
       \Delta^F=\Delta_N^F =\Delta _N \cap F \quad \mbox{ and } B_2^F =
          B_2^{N+1}\cap F .
$$

For a metric space $(X, \rho)$ and $\eps>0$ an $\eps$-net $\cal{N}$ is
a subset of $X$ such that for every $x$ in $X$ there exists $x_0\in
{\cal{N}}$ satisfying $\rho(x, x_0) \leq \eps$.

Let $k\leq d$. By $O(d)$ we denote the group of orthogonal operators
on $\R^d$ and by $G_{d, k}$ we denote the Grassmannian of
$k$-dimensional linear subspaces of $\R^d$ endowed with the distance
$$
   \rho (E, F) = \inf\{ \|U - I\| \, \mid \, U\in O(d), UE=F\} ,
$$
where $\no$ denotes the operator norm $\ell _2^d \to \ell _2^d$.

We will use the following result of Szarek (\cite{S1, S2}) on the size
of $\eps$-nets on $G_{d, k}$.

\begin{theorem}
\label{net}
  Let $k\leq d$ and $\eps \in (0, 1)$. There exists an $\eps$-net on
  $G_{d, k}$ with respect to $\rho(\cdot)$ of cardinality not
  exceeding $(C/\eps)^{C d k}$, where $C$ is an absolute positive constant.
\end{theorem}

\bigskip

Volume estimates play an important role in the theory.
Let us recall the following fundamental result (\cite{BF, CP, G1}).

\begin{fact}\label{volarg}
  \label{fact:cp}
  Let $M \ge 2d$ be integers. For
arbitrary vectors  $x_1, \ldots, x_{M} \in S^{d-1}$
 the volume of the absolute convex hull satisfies
\[
  \vol (\absconv \{x_1, \ldots, x_{M}\} )
\le   \left(C \frac{\sqrt{\ln (M/d)}}{d}\right)^d,
\]
where  $C$ is a positive absolute constant.
\end{fact}

\bigskip

The proof of existence of convex bodies that are poorly approximated
by projections of sections of a simplex uses 
a modification of bodies introduced by 
Gluskin in \cite{G2}. 
This probabilistic construction and its further versions
became the main source of counterexamples in
asymptotic geometric analysis \cite{MT}.  
However, most polytopes described in the literature have the number of
random vertices $M$  proportional to $d$, while we  want  $M$ 
to be  arbitrary satisfying  $2d \le M \le e^d$. 
To keep this paper self-contained we show
an existence with a direct argument.

Let $d \ge 1$ and $2d\leq M \leq e^d$ be integers. Set
\[
 \ell=\left \lceil \log_5 (M/d) \right \rceil,
\]
and let $\{1, \ldots, d \}=\bigcup_{k=1}^{\lceil d/\ell \rceil} I_k$
be the decomposition of $\{1, \ldots, d\}$ into the disjoint union of
consecutive intervals, with  each interval,
except possibly the last one, 
 consisting of $\ell$ numbers. 
For each $1\le k\le {\lceil d/\ell \rceil}$ choose a $(1/2)$-net $\NN_k
\subset S^{d-1} \cap \R^{I_k}$ of cardinality at most $5^\ell$. (It is
well known that such a net exists, cf. Lemma~\ref{sizenet} below;
moreover, one can show that such a net can be taken symmetric about
the origin.)

 Recall that $\P$ is the rotation invariant
probability measure on the Euclidean unit sphere $S^{d-1}$.  (We may
also denote this probability space by $(\Omega, \P)$.)  Let $X$ be a
random vector
uniformly distributed on $S^{d-1}$, and let $X_1, \ldots, X_M$ be
independent copies of $X$.  Then 
we define Gluskin's polytope $\ptop \subset \R^d$ by
\begin{equation}
  \label{eq:glu_poly}
\ptop = \absconv \left \{\bigcup_{i=1}^d \{ e_i \} \cup
  \bigcup_{k=1}^{\lceil d/\ell \rceil} \NN_k \cup  
   \bigcup_{j=1}^M \{X_j\} \right \}.
\end{equation}
To emphasize the number of random vertices we will denote $\ptop$ by
$\ptop_M$.  Since $\NN_k$ is symmetric, $2d \leq M$, and by the choice
of $\ell$, we observe that  
$\ptop_M$ has less than or equal to $4 M$ vertices. Therefore, by Fact
\ref{fact:cp}, 
\begin{equation}
  \label{eq:cp}
   \vol (\ptop_M ) \le   \left(C \frac{\sqrt{\ln (M/d)}}{d}\right)^d. 
\end{equation}
 This definition of Gluskin's polytopes differs from the original one
 in \cite{G2} by the inclusion of  
the nets $\NN_k$. 
This  guarantees that the polytope $\ptop_M$ contains a ball of 
an appropriate radius, which is necessary for the construction
below. Let $x=(x_1, \ldots, x_d) \in \R^{d}$.  
Since $\NN_k$ is a $(1/2)$-net in $S^{d-1}\cap \R^{I_k}$, 
we have
$(1/2) B_2^{I_k} \subset  
\conv (\NN_k) \subset \ptop_M$. Therefore,
$$
  \|x\|_{V_M} = \left\| \sum_{k=1}^{\lceil d/\ell \rceil} 
   \sum_{j \in I_k} x_j e_j  \right\| _{V_M} \leq   \sum_{k=1}^{\lceil d/\ell \rceil} \left\|
   \sum_{j \in I_k} x_j e_j  \right\| _{V_M}
$$
$$
  \leq  2 \sum_{k=1}^{\lceil d/\ell \rceil} \left|\sum_{j \in I_k} x_j e_j  \right| \leq 
  2 \, \sqrt{\lceil d/\ell \rceil} 
   \, \left( \sum_{k=1}^{\lceil d/\ell \rceil} \left | \sum_{j \in I_k} x_j e_j \right |^2 \right)^{1/2} 
  \leq   4 \, \sqrt{ \frac{d}{\ln (M/d)} } \ \, |x|,
$$
 which means that
\begin{equation}
  \label{eq:ball inside}
   B_2^d \subset
     4 \sqrt{\frac{d}{\ln (M/d)}}\  \ptop_M.
\end{equation}

\medskip

Having two independent Gluskin's polytopes $\ptop_M'$ and $\ptop_M''$
in $\R^d$ we will represent them on the product space 
$S^{d-1} \times S^{d-1}$ with the product probability $\P\otimes \P$.
The next theorem shows that with high probability two Gluskin polytopes 
are far apart in the Banach-Mazur distance.
The proof of this Theorem will be presented in Section \ref{appendix}.

\medskip

\begin{theorem}
  \label{mt}
  There exists a (small)  constant 
$a >0$  such that for all integers $2d\leq M \leq e^d$
  the subset of
  pairs $(\ptop_M', \ptop_M'')$ of two independent Gluskin's polytopes
  in $\R^d$ satisfies
\begin{equation}
  \label{eq:two_glu}
 \P\otimes \P\left( \left\{(\ptop_M', \ptop_M'') \mid  d(\ptop_M',
     \ptop_M'')\leq  \frac{a\, d}{\ln (M/d)} 
\right\}\right) \leq   2e^{- d M} .
  \end{equation}
\end{theorem}

\medskip

\begin{cor}   
\label{glus} 
Let $2d\leq M \leq e^d$.  
Let $K\subset \R^d$ be a convex body. 
Then Gluskin's polytopes $\ptop_M$
   in $\R^d$ with $M$ random vertices satisfy
$$
  \P\left( \left\{ \ptop_M \mid d(\ptop_M, K )\leq
      C \sqrt{\frac{ d}{\ln (\frac{M}{d})}}  \right\}\right)
         \leq \sqrt 2 e^{- d M/2},
$$
where $C >0$ is an absolute constant.
\end{cor}

\pr
Let $\ptop_M$, $\ptop_M'$, and $\ptop_M''$ be independent Gluskin's
polytopes in $\R^d$ with  $M$ random vertices. 
By Theorem~\ref{mt} and submultiplicativity of the Banach-Mazur distance, 
for every convex body $K$ we have 
\begin{align*}
 2 e^{- d M} 
 &\geq \P\otimes\P\left( \left\{ (\ptop_M', \ptop_M'') \mid   d(\ptop_M', \ptop_M'')
     \leq  \frac{a\, d}{\ln (\frac{M}{d})} \right\}\right)\\
  &\geq \P\otimes\P\left( \left\{ (\ptop_M', \ptop_M'') \mid  d(\ptop_M', K) d(K, \ptop_M'')\leq
      \frac{a\, d}{\ln (\frac{M}{d})} \right\}\right) \\
  &\geq \P\otimes\P\left(\left\{  (\ptop_M', \ptop_M'') \mid 
  \max\{  d(\ptop_M', K), d(K, \ptop_M'' )\} \le
  \sqrt{ \frac{a\, d}{\ln (\frac{M}{d})}}   
             \right\}\right) \\
 & = 
 \left(\P\left( \left\{ \ptop_M \mid d(\ptop_M, K) \le
          \sqrt{ \frac{a d}{\ln (\frac{M}{d})}} \right\}\right) \right)^2,
\end{align*}
which implies the result.
\qed

\section{Proof of the main result}
\label{main}

We start with the following lemma, which shows that it is enough to consider
only special sections of the cone $S$.

\begin{lemma}\label{sect}
  Let $m\leq N$ and let $F\subset \R^{N+1}$ be an affine subspace such
  that $\Delta_N^F$ is an $m$-dimensional body. 
Then there exists a
  linear subspace $L\subset \R^{N+1}$ such that 
   $\Delta _N^F$  has a position 
$K$ inside $L$ of the form
$$
    K = \{x\in \R^{N+1} \, \mid \,  x \in L\, \, \,  \mbox{ and }\, \,
    \, -1\leq x_i\leq m 
     \mbox{ for all }  i\leq N+1 \}.
$$
In particular,
$$
      B_2^L \subset K \subset m^{3/2} B_2^L .
$$
\end{lemma}

\pr
By (\ref{center}) there exists $a= \{a_i\} _{i=1}^{N+1}\in \Delta_N^F \subset S$ such that
\begin{equation} \label{minus}
    -(\Delta_N^F -a) \subset m (\Delta_N^F-a).
\end{equation} 
Clearly $a_i\geq 0$ for all $i\leq N+1$.
Without loss of generality we can assume that $a_i>0$ for all $i$. Indeed, note that
$a$ is in the relative interior of $\Delta^F_N$. Thus, if for some
$j>0$, $a_j=0$ then
$$
   \Delta _N^F \subset H_j:= \{ x\in \R^{N+1} \, \mid \, x_j=0\}.
$$
Therefore $\Delta _N^F$ is in fact a corresponding $m$-dimensional section of the
$(N-1)$-dimensional simplex
$$
      \Delta _{N-1} = S\cap H \cap H_j
$$
and we can apply the proof below for this section (or just to take the operator $D$ below
with zero $j$-th row).

Consider the diagonal operator $D$ with $1/a_i$'s on the main diagonal.
Denote
$$
  b=D a= \sum _{i=1}^{N+1} e_i \quad \mbox{ and } \quad K :=
  D (\Delta _N^F - a) = D \Delta _N^F - b.
$$
Then
$$
  D \Delta _N^F = D(S\cap H\cap F) = S \cap D(H\cap F).
  $$
Therefore, denoting $L:=D(H\cap F) -b$, we obtain
$$
 K= \{x\in \R^{N+1} \, \mid \, -1\leq x_i \mbox{ and } x \in L\} . 
$$
By (\ref{minus}) we observe that $-K \subset mK$, hence
$$
    K = \{x\in \R^{N+1} \, \mid \,x \in L \, \, \,  \mbox{ and }\, \, \, -1\leq x_i\leq m
     \mbox{ for all }  i\leq N+1 \}.
$$
This implies
$$
    B_2^{N+1} \cap L \subset K \subset m^{3/2} B_2^{N+1} \cap L.
$$
\qed

\begin{lemma}\label{dsect}
Let $\eps \in (0,1)$ and $m\leq N$. 
For $j=1,2$ let $L_j$ be an
$m$-dimensional linear subspace of $\R^{N+1}$  and put
$$
    K_j :=  \{x\in \R^{N+1} \, \mid \,x \in L_j \, \, \,  
\mbox{ and }\, \, \, -1\leq x_i\leq m
     \mbox{ for all }  i\leq N+1\}.
$$
Assume $\rho (L_1, L_2)\leq \eps$. 
Then  
$$
    d(K_1, K_2) \leq (1+\eps m^{3/2})^2.
$$
\end{lemma}

\pr By the definition there exists an orthogonal operator $U$ such
that $U L_1 = L_2$ and $\|U - I\|\leq \eps$. Therefore for every
$x=\{x_i\}_i \in K_1$ we have $|Ux-x|\leq \eps |x|\leq \eps m^{3/2}$,
hence $|(Ux-x)_i|\leq \eps m^{3/2}$ for every $i\leq N+1$. Thus, for
every $i$ we have
$$
   (Ux)_i = x_i + (Ux-x)_i\geq -(1+\eps m^{3/2})
$$
and
$$
   (Ux)_i = x_i + (Ux-x)_i\leq m  +  \eps m^{3/2} .
$$
Therefore, $UK_1 \subset (1+\eps m^{3/2}) K_2$. Similarly, $U^{-1} K_2
\subset (1+\eps m^{3/2}) K_1$, 
which implies the result.
\qed

\begin{lemma}\label{dproj}
Let $\eps \in (0,1)$, $n\leq m\leq N$, 
$L$ be an $m$-dimensional
linear subspace of $\R^{N+1}$ and 
$$
    K = \{x\in \R^{N+1} \, \mid \, x \in L \, \, \,  
          \mbox{ and }\, \, \,  -1\leq x_i\leq m
     \mbox{ for all }  i\leq N+1\}.
$$
Let $F_1$ and $F_2$ be $n$-dimensional linear subspaces of $\R^{N+1}$
and $P_1$ and $P_2$ be the orthogonal projections on $F_1$ and $F_2$,
respectively. Assume $\rho (F_1, F_2)\leq \eps$. Then 
$$
    d(P_1 K, P_2 K) \leq (1+\eps m^{3/2})^2 .
$$
\end{lemma}

\pr
By the definition there exists an orthogonal operator $U$ such that $U F_1 = F_2$
and $\|U - I\|\leq \eps$. Then $UP_1=P_2U$ and therefore for every $x\in K$ we have
$$
     UP_1 x = P_2 Ux = P_2 x + P_2 (U-I) x \in P_2 K + P_2 (U-I) m^{3/2}
     B_2^{N+1}\cap L .
$$
Since $B_2^{N+1}\cap L\subset K$, we obtain 
$$
   UP_1 x \in  (1+\eps m^{3/2}) P_2 K.  
$$
Similarly,
$$
     U^{-1} P_2 x    \subset (1+\eps m^{3/2}) P_1 K,
$$
which implies the result.
\qed

\medskip

We are now ready  to prove our main theorem.

\smallskip

\noindent
{\bf Proof of Theorem~\ref{mainth}:  \ }
In this proof $ C_1, C_2, C_3$ are absolute constants greater then
one.   Without loss of generality we assume 
 that $2\leq n\leq N\leq e^{cn}$, where $c$ is an absolute positive constant, which 
will be specified later (if $n=1$ or $N\geq e^{cn}$ the conclusion of the theorem 
is immediate).

For any $k\leq N$ and $\eps\in (0, 1)$, by ${\cal{A}}_k$ we denote an
$\eps$-net on the Grassmanian $G_{N+1, k}$ of cardinality 
$$
     \left| {\cal{A}}_k \r| \leq \left(C_1/\eps\r) ^{C_1 N k} .
$$
(The existence of such a  net follows from  Lemma~\ref{net}. Note that
we suppress  the dependence of the net on $\eps$.)

In the first part of the argument fix an integer $m$ such that $n\leq
m \leq N$ and fix $\eps \in (0,1)$. Put
$$
    K_m = \{x\in \R^{N+1} \, \mid \, -1\leq x_i\leq m  \mbox{ for all
    }  i\leq N+1\} . 
$$
Let    $2n \le M \le e^n$.
We apply Corollary~\ref{glus} with $d=n$
and the body $K = P_{E_0} (K_m\cap L_0)$, 
for arbitrary $L_0\in {\cal{A}}_m $  and $ E_0\in {\cal{A}}_n $.
By  the union bound we obtain that for $n$-dimensional 
 Gluskin's polytopes $\ptop_M$ one has
$$
  \P \left(\left\{ \forall L_0\in {\cal{A}}_m \, \forall E_0\in
      {\cal{A}}_n \, \, \, \, \, 
  d(\ptop_M, P_{E_0} (K_m\cap L_0))\leq 
 C_2 \sqrt{\frac{ n}{\ln (\frac{M}{d})}} 
  \right\}\right) 
$$
$$
  \leq \sqrt{2}\, 
  \left(C_1/\eps\r) ^{C_1 N m +C_1 N n}  \exp(- M n/2) \leq
  \sqrt{2}\, \exp(-M n/2  + 2 C_1 N m \ln(C_1/\eps)).
$$
Therefore whenever $M$ satisfies 
\begin{equation} 
\label{cond}
 M\geq  8 C_1 N m \ln(C_1/\eps)/  n , 
\end{equation}
then 
$$
    \P  \left(\left\{ \forall L_0 \in {\cal{A}}_m  \, \forall E_0 \in
      {\cal{A}}_n  \, \, \, \, \, \,  d(\ptop_M, P_{E_0} (K_m\cap L_0))\leq 
    C_2 \sqrt{\frac{ n} {\ln (\frac{M}{n})}}\, \right\}\right) 
$$
\begin{equation}\label{cond_prob}
   \leq  \sqrt{2} \, \exp(- M n/4 )\leq \exp(-M n/6).
\end{equation}


Therefore  taking  $M$  satisfying $2n \le M\le e^n$  and (\ref{cond})
(if such an $M$ exists), 
this implies the result for 
Gluskin's polytopes $\ptop_M$ and for 
every $n$-dimensional projection of an
$m$-dimen\-sio\-nal section of an $N$-dimensional simplex, with high
probability.  (Note that $m$ is fixed in this argument.)  Indeed, let
$F$ be any affine subspace of $\R^{N+1}$, such that $\Delta _N^F$ is
$m$-dimensional.  Let $L=L(F)$ be an $m$-dimensional linear subspace
and $K=K(F)= K_m\cap L$ be the position of $\Delta _N^F$ provided by
Lemma~\ref{sect}. Let $P$ be any orthogonal projection such that $PK$
is $n$-dimensional and let $E$ be the range of $P$. Let $L_0\in
{\cal{A}}_m$ and $E_0\in {\cal{A}}_n$ be such that $\rho(L, L_0) \leq
\eps$ and $\rho(E, E_0) \leq \eps$. Then by Lemmas~\ref{dsect} and
\ref{dproj} we get
$$
  d(P K, P_{E_0} (K_m\cap L_0)) \leq  d(P K, P_{E_0} K)  d(P_{E_0} K,
  P_{E_0} (K_m\cap L_0)) 
  \leq  (1+\eps m^{3/2})^4, 
$$
where in the last estimate we used the obvious inequality $d(P_{E_0}
K_1, P_{E_0} K_2) \le d(K_1,K_2)$  
valid for all convex bodies $K_1, K_2 \subset \R^N$ of dimension $m \leq N$.
Therefore, taking $\eps = m^{-3/2}$ we obtain that 
$$
 d(\ptop_M,  P_{E_0} (K_m\cap L_0))   \le 2^{4} d(\ptop_M, PK).
$$
Combining this with (\ref{cond_prob}), we 
obtain the probability  estimate for 
$$
  \P \left(\left\{ \mbox{for every }\, F, L, K, P \, \mbox{ as above: } \, \, 
   d(\ptop_M, PK)\leq 2^{-4} 
 C_2 \sqrt{\frac{ n}{\ln (\frac{M}{n})}} 
\right\}\right) .
$$
More precisely we showed that  for any $n\le m \le N$
whenever $M$ satisfies $2n \le M \le e^n$ and  (\ref{cond}) 
with $\eps = m^{-3/2}$, 
then the latter probability is less than or equal to 
$ \exp(-   M n/6 )$.  
In particular, let  
$$
M=\lceil 8 C_1 N^2 \ln(C_1 N^{3/2})/ n  \rceil ,  
$$
so that (\ref{cond}) is satisfied 
with $\eps = m^{-3/2}$. 
Additionally we can find a universal constant $0 < c <1$
such that the condition $N \le e^{c\,n}$ implies 
$M\le e^n$.

Then for some absolute constant $C_3$,
$$
  \P \Bigl(\Bigl\{ \mbox{for every }\, K, P, \
d(\ptop_M, PK)\leq C_3
 \sqrt{  \frac{ n }{ \ln \frac{2 N \ln (2N)}{n} }}
\Bigr\}\Bigr)
  \leq   \exp(-  N^2  \ln(2 N)).
$$
(Here $K$ and $P$ are as above, in particular, the dimension of a
section $K$ is equal to $m$.)

To obtain the full result for any $n\le N$,  for any
$n$-dimensional projection of an arbitrary dimensional section of an
$N$-dimensional simplex we apply the above discussion for an arbitrary
$m$ representing the dimension of a section (so $n\le m\le N$). Note
that the choice of $M$ does not depend on $m$, so we are working in
the same probability space for all $m$, leading to the same class of
Gluskin's polytopes $\ptop_M$.  Taking the union bound over all
integers $n\le m\le N$ concludes the proof.  
\qed

\medskip

\noindent
{\bf Remarks. 1.} In fact, taking $M=\lceil 8 C_1 N m \ln(C_1 m^{3/2})/ n  \rceil$
in our proof, we observe  that for $n\leq m\leq N$
there exists an $n$-dimensional convex body $B$ such that for every
convex body $K$ obtained as an $n$-dimensional projection of an
$m$-dimensional section of an $N$-dimensional simplex one has
$$
    d(B, K) \geq  c \sqrt{  \frac{ n }{ \ln \frac{2 N m \ln (2 m)}{n^2}  }} . 
$$
Moreover, our construction is random -- we use 
Gluskin's polytopes --
and we obtain the result with high probability -- the estimate above
holds with probability larger than $1-\exp(- N m \ln (2m))$.

\noindent
{\bf 2. }  If we restrict ourselves to just one operation --
projection -- then we have almost the same lower bound using the
Euclidean ball. Namely, for every $n$-dimensional projection $P$ one
has
 $$
    d(B_2^n, P\Delta _N) \geq c  \sqrt{  \frac{ n }{ \ln \frac{2 N }{n}  }} , 
$$
which follows from volume estimates (see Fact~\ref{volarg}) as  mentioned
in the introduction.

\noindent
{\bf 3. }
Also note that, although an $N$-dimensional simplex clearly has $\lceil N/2 \rceil$-dimensional
symmetric projection, a ``random" projection is very far from being
symmetric. It was shown in  Theorem~5.1 of \cite{LT} that for a ``random"
$n$-dimensional projection $P$ and every centrally symmetric convex body $B$
one has
$$
    d(B, P\Delta _N) \geq c  \sqrt{  \frac{ n }{ \ln N   }} .
$$

\bigskip

\section{Proof of Theorem~\ref{mt}}
\label{appendix}

The proof of the theorem is standard and follows the road-map of
\cite{G2}.  The main difference from \cite{G2} is the modification of
the definition of a Gluskin polytope \eqref{eq:glu_poly}.  Adding the
nets $\mathcal{N}_k$ to the vertex set of $\ptop_M$ allowed to
guarantee the inclusion \eqref{eq:ball inside} without significantly
increasing the number of vertices. 
(Of course if the number of vertices is  proportional then 
\eqref{eq:ball inside} is automatically satisfied.)

Recall that the underlying probability space is the product space 
$\Omega'\times \Omega'' = S^{d-1} \times S^{d-1}$ with the product
probability $\P\otimes \P$.  Our first aim in the proof is to prove two
estimates similar to \eqref{eq:two_glu}: one is for probability
on $\Omega'$, with $\omega'' \in \Omega''$ fixed, and in the other one
the roles of $\Omega'$ and $\Omega''$ are interchanged.  This is
proved in Lemma~\ref{half-glusk} below.  Then the full
Theorem~\ref{mt} follows by considerations based on Fubini's theorem.

Throughout most of this section, until the final proof of the theorem,
we fix an arbitrary $\omega''\in \Omega''$ and the corresponding
Gluskin's polytope $W_M = \ptop_M''(\omega'')$.

\medskip

For any $\tau >0$ and any operator $T:\R^d \to \R^d$ with $\det T =1$ 
consider the event
%
%
\begin{equation}  \label{eq:A_T}
A(\tau, W_M, T) = \big\{ \ptop_M: \ \|T: \ptop_M \to W_M\|\le \tau \big\} 
= \big\{ \ptop_M: \ T\ptop_M \subset \tau  W_M \big\}.  
\end{equation}

First we estimate the probability of this event.

\begin{lemma}
  \label{single-op}
  One has 
$$
  \P\left(A(\tau, W_M, T)\right) \le
  \bigl(C \tau \sqrt{\ln (M/d)/d}\, \bigr)^{dM}, 
$$
where $C$ is a positive absolute constant.
\end{lemma}

To prove this lemma we need the following well-known simple fact, which 
 can be found in many places, for example in \cite{Tom},  (38.4). 
We outline the proof for the reader's convenience.

\begin{fact}
  \label{measure_sphere}
Let $K \subset \R^d$ be a convex body with $0$ in its interior.
Let  $X$  be a random vector uniformly distributed on the sphere
$S^{d-1}$. Then
$$
  \P(\{X\in K\}) \le {{\vol\, (K)}/{\vol\, (B_2^d)}}.
$$
\end{fact}

\medskip

\pr
Obviously  we have
$ \P (\{X \in K\}) =  \vol(L)/\vol (B_2^d)$
where $L = \left\{x\in B_2^d \mid x/|x| \in K\cap S^{d-1}\right\} $.
On the other hand,
$ L \subset \conv (K\cap S^{d-1}) \subset  K$, 
which yields the required estimate for volumes.
\qed

\medskip

We use
a convenient shortcut for norms of linear operators:
for two convex bodies $K_1, K_2 \subset \R^d$ and for $\lambda >0$  
the statement $\|T:K_1 \to K_2\|\le \lambda$ is equivalent to
$T(K_1) \subset \lambda K_2$ and is 
equivalent to  $\|T:K_1 \to \lambda K_2\|\le 1$.

\bigskip

{\noindent {\bf Proof of Lemma~\ref{single-op}: }} Since $V_M$ contains the
 vectors $X_j, \ j\leq M$,
 the condition $T(V_M) \subset \tau W_M$
implies  that $TX_j\in \tau W_M$ for all $j\leq M$.
Therefore
\begin{eqnarray*}
\P\left(A(\tau, W_M, T)\right)
&\le &
\P\left(\left\{TX_j \in \tau  W_M
\quad {\rm for\ } 1 \le j\le M \right\}\right) \\ 
 &=&
\Big(\P\left(\left\{X \in  \tau \, T^{-1}\, W_M\right\}\right)\Big)^{M} 
\end{eqnarray*}
(cf. Lemma 38.3 in \cite{Tom} and  Lemma 4 in \cite{MT}). 
By Fact~\ref{measure_sphere} and using $\det T^{-1}=1$ and (\ref{eq:cp}) for $W_M$,  
we obtain 
$$
 \P\left(A(\tau, W_M, T)\right) \leq 
 \left(\frac{\vol (\tau   W_M)}{\vol (B_2^d)}\right)^{M}
= \tau ^{dM}  \left(\frac{\vol (W_M)}{\vol (B_2^d)}\right)^{M}
 \le \Big(C \tau \sqrt{\frac{\ln(M/d)} {d}}\,\Big)^{d M} ,
$$
which completes the proof.
\qed

\bigskip

In the next step  we discretize  certain sets
of operators acting on $\R^d$ (see 
Lemma 38 in \cite{Tom} and Lemma 7 in  \cite{MT}). We need more  notation.
Set
$$
B_{op}^d = \{T: \R^d \to \R^d \mid
 \|T: \ell_2^d\to \ell_2^d\| \le 1\},
$$
and for a convex body $K \subset \R^d$,
$$
B_{op,K}^d = \{T: \R^d \to \R^d  \mid
\|T: B_1^d\to K\| \le 1  \}.
$$
Note that the norm for which $B_{op,K}^d$ is the unit ball is equal to
the $\ell_\infty$-direct-sum of $d$ norms $\|\cdot\|_K$ determined by $K$.

For the reader's convenience we recall that identifying the set of
operators with $\R^{d^2}$ we have
\begin{equation}
  \label{eq:volumes_sets-oper}
\vol (B_{op,K}^d) = (\vol(K))^d
\qquad {\rm and } \qquad
\vol(B_{op}^d) \ge (c/ \sqrt d)^{d^2},
\end{equation}
where $c$ is a positive absolute constant.

We also will use the following fact on cardinality of $\eps$-nets. 
Recall that the smallest cardinality of a $1$-net of a set $K_1$ 
in the metric of defined by a convex body $K_2$ is denoted by 
$N(K_1,  K_2)$, hence the smallest cardinality of an $\ep$-net 
is $N(K_1, \ep K_2)$. 
The following lemma follows by the standard volumetric argument 
(in such a formulation it is Lemma~6 from \cite{MT}).

\begin{lemma} \label{sizenet}
Let $\ep >0$. Let $K_1, K_2 \subset \R^n$ be two symmetric convex bodies such that  
$K_1\subset K_2$. Then every subset $K' \subset K_2$ admits an $\ep$-net $\mathcal{N} \subset K'$ in
the metric of $K_1$ with $|\mathcal{N}|\le (1+ 2/\ep)^{n} \left(\vol( K_2)/ \vol (K_1)\right)$. 
\end{lemma}

We use this lemma to control the cardinality of an $\ep$-net in $B_{op. \eta {W}} ^d$ in the operator norm.

\begin{lemma}
  \label{lemma_net}
Let $\xi >0$ and let $W \subset \R^d$ be a convex symmetric body such that $B_2^d \subset \xi W$.
Let $\eta, \ep  >0$.
Every subset $K'$ of
$B_{op. \eta {W}} ^d$ admits an $\ep$-net
$\mathcal{N} $ in $K'$ in the operator norm on
$\ell_2^d$ with cardinality
\begin{equation}
  \label{eq:net_estimate}
  |\mathcal{N}|\le \left(\frac{\xi}{\eta}+\frac{2}{\ep}\right)^{d^2}\,
\left(C\, \eta \sqrt{d} \cdot \vol^{1/d}(W)\right)^{d^2},
\end{equation}
where $C$ is an absolute positive constant.
\end{lemma}

\pr
We will use Lemma~\ref{sizenet} with 
$\lambda = \xi/\eta$, $K_1= (1/\lambda) B_{op}^d$ and $K_2 = B_{op, \eta W}^d$.
The assumption $B_2^d \subset \xi W$ yields
$(1/\lambda) B_{op}^d \subset B_{op, \eta W}^d$. 
Thus, by (\ref{eq:volumes_sets-oper}),
\begin{align*}
N\Big(K',  \ep B_{op}^d \Big)=
N\left(K', \lambda \ep \left(\left(1/\lambda\right) B_{op}^d\right)\right)
&\le
\left(1+ \frac{2}{\ep \lambda} \right)^{d^2} \ \frac{\vol K_2}{\vol K_1} \\
&\le \left(\frac{2}{\ep} +
\frac{\xi}{\eta}\right)^{d^2} \
 \left(C\, \eta \sqrt{d} \cdot \vol^{1/d}(W)\right)^{d^2},
\end{align*}
with an absolute positive constant $C$. 
\qed

\medskip

We need one more lemma, which estimates the probability of the following event 
\begin{equation}
  \label{eq:tilde_A_no}
\widetilde A(\eta, W_M) =
  \big\{ \ptop_M : \ \exists\   S:\R^d \to \R^d, \, \det S =1,\
   {\rm s.t.\ }
 \|S: \ptop_M \to W_M\| \le \eta \big \} ,
\end{equation}
where $\eta$ is a positive parameter.

\bigskip
\begin{lemma}
  \label{half-glusk} 
Let $d\leq M \leq e^d$. There exists a positive constant $a_1>0$ such that for 
$\eta = a_1 \sqrt{d/\ln (M/d)}$ one has 
\[
  \P \Big( \widetilde A(\eta, W_M)  \Big) \le e^{-d M}. 
\] 
\end{lemma}

\pr
Denote for shortness $\xi = 4 \sqrt{\frac{d}{\ln (M/d)}}$.
Fix an arbitrary $0 < \ep \le 1$. By $K'$ denote the set
of all operators $T \in B_{op, \eta{W_M}}^d$ with
$\det T =1$.
Let $\mathcal{N}$ be an  $\ep$-net for $K'$
with respect to  the metric given by  $B_{op}^d$
and satisfying (\ref{eq:net_estimate}) with $W=W_M$.

We  first  show that
\begin{equation}
  \label{eq:tilde_A_incl}
  \widetilde A (\eta, W_M) \subset \bigcup_{T\in \mathcal{N}}
A(\tau, W_M, T),
\end{equation}
where $\tau = \eta + \ep \xi$.

Pick  $\omega \in \widetilde A(\eta, W_M)$,
and let $S$ be an operator with $\det S =1$ such that
$ \|S: \ptop_M (\omega)\to W_M\| \le \eta $.
Since $\ptop_M \supset B_1^d$, we have
$\|S: B_1^d \to \eta W_M\|\le 1$, which means
$ S \in B_{op, \eta W}^d$.

\smallskip

Since $\det S =1$, then $S$ belongs to $K'$. By the definition of
$ \mathcal{N}$, we can find
$T \in \mathcal{N}$ satisfying  $\|T-S: \ell_2^d \to \ell_2^d \|\le \ep$.
Since $\ptop_M \subset B_2^d$ and by \eqref{eq:ball inside}, we get 
$$
 (T-S)(\ptop_M)\subset \ep B_2^d \subset \ep \xi W_M.
$$
Equivalently,
$\|T-S: \ptop_M \to W_M\|\le \ep \xi$.
By the triangle inequality,
$$
\|T: \ptop_M \to W_M\|\le \|T-S: \ptop_M \to W_M\|
+ \|S: \ptop_M \to W_M\|\le \ep \xi + \eta = \tau.
$$
This means that  $\omega \in A(\tau, W_M, T)$ for every 
$T \in  \mathcal{N}$ and ends the proof
of (\ref{eq:tilde_A_incl}).

\smallskip

By  the union bound and Lemma \ref{single-op}
\[
\P\big(\widetilde A (\eta, W_M) \big) \le |\mathcal{N}|
  \left( C \tau \sqrt{\frac{\ln (M/d)}{d}}\, \right)^{d M}.
\]
Combining this with (\ref{eq:net_estimate}), (\ref{eq:cp}) for $W_M$, 
and the definitions of $\xi$ and $\eta$ we observe  that
\begin{align*}
  \P\big(\widetilde A (\eta, W_M) \big) 
  &\le
    \left(\frac{2}{\ep} + \frac{\xi}{\eta}\right)^{d^2}\,
   \left(C_1\, \eta \sqrt{d} \cdot \vol^{1/d}(W_M)\right)^{d^2} \cdot
   \left(C (\eta+\ep\xi) \sqrt{\frac{\ln (M/d)}{d}}\,
  \right)^{d M} \\
 &\le
    \left(\frac{2}{\ep} + \frac{4}{a_1}\right)^{d^2} \, (C_2 a_1)^{d^2} \cdot
    \left( C_3  (a_1+4\ep)   \right)^{d M}   , 
\end{align*}
where $C, C_1, C_2, C_3$ are absolute positive constants. 
To complete the proof it is enough to set $\ep=a_1$ and choose $a_1$ appropriately small.
\qed

\bigskip

Now we are ready to prove Theorem~\ref{mt}. 

\medskip

\noindent
{\bf Proof of Theorem~\ref{mt}: }
Let $a_1$ and $\eta$ be as in Lemma \ref{lemma_net}.
We consider various subsets of the measure spaces
$\Omega'$, $\Omega''$, and $\Omega' \times\Omega''$; we will use an
expanded notation to avoid confusion.

Denote the set that appears in
\eqref{eq:two_glu} by $D$, that is
$$
   D = \left\{(\omega',\omega'') \mid d(\ptop_M'(\omega'),
  \ptop_M''(\omega''))\leq \eta^2 \  \right\} .
$$

For any $\omega_0''\in \Omega''$ define the subset $D'_{\omega_0''}
\subset \Omega'\times\Omega''$ which depends only on the first variable
$\omega'$ with the second variable fixed $\omega'' = \omega_0''$ and
is given by
\[
D'_{\omega_0''} =
  \big\{ (\omega', \omega_0'') \mid \exists\
S \  {\rm s.t.\ }  \det S =1\
   {\rm and\ }
 \|S: \ptop_M'(\omega')\to \ptop_M''(\omega_0'')\| \le \eta \big \} .
\]

Similarly, for
any $\omega_0'\in \Omega'$ define the subset $D''_{\omega_0'}$ by
\[
D''_{\omega_0'} =
  \big\{  (\omega_0', \omega'') \mid \exists\
R \  {\rm s.t.\ }  \det R =1\
   {\rm and\ }
 \|R: \ptop_M''(\omega'')\to \ptop_M'(\omega_0')\| \le \eta \big \} .
\]
Note that both definitions closely follow the model of
(\ref{eq:tilde_A_no}) in that the norm of operators is considered from
a random polytope to a fixed polytope.

The following
inclusion can be easily checked
\[
     D  \subset \bigcup_{\omega_0'\in \Omega'} D''_{\omega_0'} \ \cup\
      \bigcup_{\omega_0''\in \Omega''} D'_{\omega_0''} .
\]
Indeed, if $d(\ptop_M'(\omega_0'), \ptop_M''(\omega_0''))\leq \eta^2$ then there 
exists an invertible operator $S$ such that 
$$
  \|S: \ptop_M'(\omega_0')\to \ptop_M''(\omega_0'')\| \, 
  \|S^{-1}: \ptop_M''(\omega_0'')\to \ptop_M'(\omega_0')\| \le \eta ^2. 
$$
Without loss of generality we may assume that $\det S = \det S^{-1} =1$. 
Thus one of the norms in the above product is less than or equals to $\eta$, 
which means that either $(\omega_0', \omega_0'') \in D'_{\omega_0''}$
or $ (\omega_0', \omega_0'') \in D''_{\omega_0'}$.

Finally, using Lemma~\ref{half-glusk} and the Fubini theorem, we obtain 
\[
 \P \otimes \P(D) 
 \le \E_{\omega_0'} \P(D''_{\omega_0'} \mid \ \omega_0')
 +\E_{\omega_0''} \P(D'_{\omega_0''} \mid \ \omega_0'')
 \le 2 e^{-d M}.
\]
 This  completes the proof of Theorem~\ref{mt}.  \qed

\vspace{1cm}

\address

\end{document}